\documentclass[reqno,11pt]{amsart}
\usepackage[letterpaper,hmargin=1.25in,vmargin=0.9in]{geometry}
\usepackage[T1]{fontenc}
\usepackage{lmodern}
\usepackage{microtype}
\usepackage{mathtools,amssymb}
\usepackage{tikz}
\usepackage{placeins}
\usepackage[unicode,hidelinks]{hyperref}
\allowdisplaybreaks[1]
\newcommand{\R}{\mathbb{R}}
\newcommand{\Z}{\mathbb{Z}}
\DeclareMathOperator{\conv}{conv}
\DeclareMathOperator{\interior}{int}
\DeclareMathOperator{\vol}{vol}
\DeclareMathOperator{\GL}{GL}
\theoremstyle{plain}
\newtheorem{theorem}{Theorem}[section]
\newtheorem{proposition}[theorem]{Proposition}
\newtheorem{lemma}[theorem]{Lemma}
\newtheorem{corollary}[theorem]{Corollary}
\newtheorem{conjecture}{Conjecture}[section]
\theoremstyle{definition}
\newtheorem{definition}[theorem]{Definition}
\theoremstyle{remark}
\newtheorem{remark}[theorem]{Remark}
\numberwithin{equation}{section}

\title[Volume, lattice points, and successive minima]{Bounds for lattice points in terms of volume and successive minima}
\author{Jun Wang}
\address{School of Mathematical Sciences \& Institute of Mathematics and Interdisciplinary Sciences, Tianjin Normal University, Tianjin 300387, China}
\email{junwang@tjnu.edu.cn}

\author{Fei Xue}
\address{School of Mathematical Sciences, Nanjing Normal University, Nanjing 210046, China}
\email{05429@njnu.edu.cn}
\date{}
\hypersetup{pdftitle={Bounds for lattice points in terms of volume and successive minima},pdfauthor={Jun Wang and Fei Xue},pdfsubject={Lattice points, volume, and successive minima},pdfkeywords={lattice points, successive minima, convex bodies, anti-blocking bodies}}

\begin{document}
\begin{abstract}
We prove, in every dimension, the upper bound conjectured by Betke, Henk and Wills for the lattice-point enumerator in terms of volume and successive minima. For every convex body $K\subset\R^n$ and full-rank lattice $\Lambda$, we establish
$$\#(K\cap\Lambda)\le\frac{\vol_n(K)}{\det\Lambda} \prod_{i=1}^n\left(1+\frac{i\lambda_i}{2}\right),$$
where $\lambda_i$ are the successive minima of $(K-K)/2$ with respect to $\Lambda$. Equality holds precisely for lattice translates of positive integral dilates of unimodular simplices. We also prove
$$\#(\interior K\cap\Lambda)\ge\frac{\vol_n(K)}{\det\Lambda} \left(1-\frac{n}{4}\sum_{i=1}^{n-1}\lambda_i-\frac{n}{2}\lambda_n\right),$$
which implies the lower estimate with factor $1-\frac{1}{2}\sum_{i=1}^n i\lambda_i$. In dimension two we prove the lower bound with factor $(1-\lambda_1/2)(1-\lambda_2)$ when $\lambda_2\le1$. We also observe that a conjectured product lower bound for origin-symmetric bodies fails in dimension three. The proofs combine geometric compression, lattice-cell counts in coordinate sections, and an elementary coefficient comparison.
\end{abstract}
\maketitle

\section{Introduction}\label{sec:1}

Let $K\subset\R^n$ be an $n$-dimensional convex body, i.e. a compact and convex set with nonempty interior. For a measurable subset $M\subset\R^n$, we denote by $\vol_n(M)$ its volume, that is, its Lebesgue measure. Let $\Lambda\subset\R^n$ be a full-rank lattice with determinant $\det\Lambda$. For an origin-symmetric convex body $C$, the $i$-th successive minimum with respect to $\Lambda$ is defined as
$$\lambda_i(C,\Lambda)=\min\{\rho>0:\dim(\rho C\cap\Lambda)\ge i\}, \qquad 1\le i\le n.$$
For an arbitrary convex body $K$, put $K_s=\frac{1}{2}(K-K)$. We use the successive minima of the origin-symmetric body $K_s$ throughout the paper.

A classical and fundamental result in geometry of numbers called Minkowski's second theorem gives, in its form for arbitrary convex bodies \cite[(1.2)]{ref09},
\begin{equation}\label{eq:1-1}
 \frac{2^n}{n!}\le\frac{\vol_n(K)}{\det\Lambda} \prod_{i=1}^n\lambda_i(K_s,\Lambda)\le2^n.
\end{equation}
Its original formulation, dating back to 1896, concerned only origin-symmetric convex bodies, see \cite[Theorem 16]{ref15} for reference. Betke, Henk and Wills \cite{ref05} proposed a discrete analogue in 1993, replacing volume by the lattice-point enumerator $\mathrm{G}(K,\Lambda)=\#(K\cap\Lambda)$. Their conjecture was originally stated for origin-symmetric bodies and later extended to arbitrary convex bodies, see \cite[Conjecture 5.2]{ref01}. In this general form, the conjecture reads
\begin{equation}\label{eq:1-2}
 \mathrm{G}(K,\Lambda)\le\prod_{i=1}^n
 \left\lfloor\frac{2}{\lambda_i(K_s,\Lambda)}+1\right\rfloor.
\end{equation}
\pagebreak[3]
For partial results toward this conjecture, see Henk \cite{ref12}, Malikiosis \cite{ref14} and Tointon \cite{ref16}. The conjecture is still open for all $n\ge4$. It is clear that a bound of the form $\prod_i(2/\lambda_i(K_s,\Lambda)+c_i)$, where the $c_i$ are independent of $K$, yields Minkowski's upper bound~\eqref{eq:1-1} by elementary properties of the Riemann integral. In order to study the discrete analogous conjecture, Betke, Henk and Wills (\cite[Conjecture 2]{ref09}; see also \cite[Conjecture 5.3]{ref01} and \cite[Conjecture 2]{ref13}) further proposed another conjecture that relates the volume, the lattice point enumerator and the successive minima simultaneously.

\begin{conjecture}[Betke--Henk--Wills]\label{conj:1-1}
For every convex body $K\subset\R^n$ and every full-rank lattice $\Lambda$, we have
\begin{equation}\label{eq:1-3}
 \mathrm{G}(K,\Lambda)\le\frac{\vol_n(K)}{\det\Lambda}
 \prod_{i=1}^n\left(1+\frac{i\lambda_i(K_s,\Lambda)}{2}\right).
\end{equation}
If $\lambda_n(K_s,\Lambda)\le\frac{2}{n}$, then
\begin{equation}\label{eq:1-4}
 \mathrm{G}(\interior K,\Lambda)\ge\frac{\vol_n(K)}{\det\Lambda}
 \prod_{i=1}^n\left(1-\frac{i\lambda_i(K_s,\Lambda)}{2}\right).
\end{equation}
Moreover, if $K=-K$ and $\lambda_n(K_s,\Lambda)\le2$, we have
\begin{equation}\label{eq:1-5}
 \mathrm{G}(\interior K,\Lambda)\ge\frac{\vol_n(K)}{\det\Lambda}
 \prod_{i=1}^n\left(1-\frac{\lambda_i(K_s,\Lambda)}{2}\right).
\end{equation}
\end{conjecture}
Apart from yielding discrete versions of Minkowski's second theorem from~\eqref{eq:1-3}, Conjecture~\ref{conj:1-1} is of interest in its own right.
For fixed $K$ and $\Lambda$, $\mathrm{G}(tK,\Lambda)/t^n$ tends to $\vol_n(K)/\det\Lambda$ as $t\to\infty$. At a fixed scale, however, lattice points on the boundary and thin directions relative to the lattice can produce a substantial discrepancy. The conjectured upper bound~\eqref{eq:1-3} controls that discrepancy by the successive minima.

In 2022, Freyer and Lucas \cite{ref09} proved~\eqref{eq:1-3} in dimension two and obtained in every dimension the bound
$$\mathrm{G}(K,\Lambda)\le\frac{\vol_n(K)}{\det\Lambda} \prod_{i=1}^n\left(1+\frac{n\lambda_i(K_s,\Lambda)}{2}\right).$$
For $\Lambda=\Z^n$, they observed that $m\Delta_n$ attains the conjectured bound for each positive integer $m$, where $\Delta_n=\conv(\mathbf{0},\mathbf{e}_1,\ldots,\mathbf{e}_n)$ is a {\it unimodular} simplex and the $\mathbf{e}_i$ are the standard basis vectors. A simplex whose edge vectors at a lattice vertex form a basis of the lattice is called {\it unimodular}.

In this paper, we first give an affirmative answer to the upper bound~\eqref{eq:1-3} in Conjecture~\ref{conj:1-1} and completely characterize the equality cases.

\begin{theorem}\label{thm:1-1}
For every convex body $K\subset\R^n$ and every full-rank lattice $\Lambda$, one has
$$\mathrm{G}(K,\Lambda)\le\frac{\vol_n(K)}{\det\Lambda} \prod_{i=1}^n\left(1+\frac{i\lambda_i(K_s,\Lambda)}{2}\right).$$
Equality holds if and only if
$$K=\mathbf{t}+m\conv(\mathbf{0},\mathbf{b}_1,\ldots,\mathbf{b}_n),$$
where $\mathbf{t}\in\Lambda$, $m\in\Z_{>0}$, and $\mathbf{b}_1,\ldots,\mathbf{b}_n$ form a basis of $\Lambda$.
\end{theorem}

Combining Theorem~\ref{thm:1-1} with Minkowski's upper volume bound gives
$$\mathrm{G}(K,\Lambda)\le\prod_{i=1}^n \left(\frac{2}{\lambda_i(K_s,\Lambda)}+i\right).$$
But this consequence is weaker than the conjectured~\eqref{eq:1-2}.

As for the lower bound for interior lattice points enumerator $\mathrm{G}(\interior K,\Lambda)$, we prove Theorem~\ref{thm:1-2} and Corollary~\ref{cor:1-3}. In particular, we prove the planar case of~\eqref{eq:1-4} in Theorem~\ref{thm:1-4}.

\begin{theorem}\label{thm:1-2}
For every convex body $K\subset\R^n$ and every full-rank lattice $\Lambda$, one has
\begin{equation}\label{eq:1-6}
 \mathrm{G}(\interior K,\Lambda)\ge\frac{\vol_n(K)}{\det\Lambda}
 \left(1-\frac{n}{4}\sum_{i=1}^{n-1}\lambda_i(K_s,\Lambda)
 -\frac{n}{2}\lambda_n(K_s,\Lambda)\right).
\end{equation}
\end{theorem}

\begin{corollary}\label{cor:1-3}
Under the same assumptions as in Theorem~\ref{thm:1-2},
\begin{equation}\label{eq:1-7}
 \mathrm{G}(\interior K,\Lambda)\ge\frac{\vol_n(K)}{\det\Lambda}
 \left(1-\frac{1}{2}\sum_{i=1}^n i\lambda_i(K_s,\Lambda)\right).
\end{equation}
\end{corollary}

For $n\ge3$, the right-hand sides in Theorem~\ref{thm:1-2} and Corollary~\ref{cor:1-3} agree precisely when $\lambda_1(K_s,\Lambda)=\cdots=\lambda_{n-1}(K_s,\Lambda)$; otherwise Theorem~\ref{thm:1-2} is stronger. Corollary~\ref{cor:1-3} proves the linear part of the conjectured bound~\eqref{eq:1-4}. After division by $\vol_n(tK)/\det\Lambda$, the upper bound and~\eqref{eq:1-7} for $tK$ have first-order corrections of the same magnitude, namely $\frac{1}{2t}\sum_{i=1}^n i\lambda_i(K_s,\Lambda)$. Integral dilates of the standard simplex show that the coefficient $1/2$ in~\eqref{eq:1-7} cannot be decreased; see Sect.~\ref{sec:7}.

\begin{theorem}\label{thm:1-4}
Let $K\subset\R^2$ be a convex body and let $\Lambda$ be a full-rank lattice. If $\lambda_2(K_s,\Lambda)\le1$, then
\begin{equation}\label{eq:1-8}
 \mathrm{G}(\interior K,\Lambda)\ge\frac{\vol_2(K)}{\det\Lambda}
 \left(1-\frac{\lambda_1(K_s,\Lambda)}{2}\right)
 \left(1-\lambda_2(K_s,\Lambda)\right).
\end{equation}
\end{theorem}

Equality is attained by $m\conv(\mathbf{0},\mathbf{b}_1,\mathbf{b}_2)$ for every integer $m\ge2$, where $\mathbf{b}_1,\mathbf{b}_2$ form a basis of $\Lambda$. Freyer~\cite{ref08} proved the corresponding estimate for real axial right triangles.

For origin-symmetric bodies, Conjecture~\ref{conj:1-1} states the stronger lower bound~\eqref{eq:1-5}. Freyer and Lucas proved this bound for origin-symmetric lattice polygons \cite[Proposition 1.5]{ref09}. To one's surprise,~\eqref{eq:1-5} fails in dimension three, as we explain in Remark~\ref{rem:8-2} using a lattice parallelepiped.

The paper is organized as follows: in Sect.~\ref{sec:2}, we first normalize the lattice to integer lattice and recall that it is enough to prove the upper bound~\eqref{eq:1-3} and the nonsymmetric lower bounds~\eqref{eq:1-4} in Conjecture~\ref{conj:1-1} for anti-blocking bodies based on the Freyer--Lucas reduction. Sect.~\ref{sec:3} introduces diagonal compressions in coordinate planes and constructs a common fixed body that is fully shifted, has the same volume, and satisfies the required lattice-point inequalities. The transfer property of fully shifted bodies allows us to convert lattice-point counts into a weighted sum of section volumes via a cell construction in Sect.~\ref{sec:4}. Combined with the ordered-sequence inequality introduced in Sect.~\ref{sec:5}, we prove the conjectured upper bound and completely characterize the equality cases of Theorem~\ref{thm:1-1} in Sect.~\ref{sec:6}. By assigning upper cells to lower cells in codimension-one sections, we give the lower bound for interior lattice points enumerator and prove Theorem~\ref{thm:1-2} and Corollary~\ref{cor:1-3} in Sect.~\ref{sec:7}. Finally, Sect.~\ref{sec:8} proves the planar product lower bound in Theorem~\ref{thm:1-4}, and presents a three-dimensional lattice parallelepiped counterexample that disproves the stronger symmetric-body conjecture~\eqref{eq:1-5}.

\section{Reduction to anti-blocking bodies}\label{sec:2}

The goal of this section is to reduce the general problem for every convex body $K$ and full-rank lattice $\Lambda$ to the case of an arbitrary anti-blocking body $B$ and the integer lattice $\Z^n$. After this reduction, the successive minima will be given explicitly by the axial intercepts of $B$.

By an invertible linear transformation of coordinates, we may assume that $\Lambda=\Z^n$. For simplicity, we omit the lattice from $\mathrm{G}$ and $\lambda_i$ when it is $\Z^n$, and we use $[n]=\{1,\ldots,n\}$. Throughout the paper, vector inequalities are understood componentwise. That is to say, for $\mathbf{x}=(x_1,\ldots,x_n)$, $\mathbf{y}=(y_1,\ldots,y_n)\in\R^n$, we write $\mathbf{y}\le\mathbf{x}$ if $y_i\le x_i$ for all $i\in[n]$.

Then we give the definition of {\it anti-blocking body}. A convex body $B\subset\R_{\ge0}^n$ is called an anti-blocking body if for every $\mathbf{x}\in B$ and $\mathbf{y}$ such that $\mathbf{0}\le\mathbf{y}\le\mathbf{x}$, we have $\mathbf{y}\in B$. Anti-blocking bodies were firstly named by Fulkerson \cite{ref10,ref11} in the context of combinatorial optimization and later studied in depth in \cite{ref02}. The projection of $B$ onto any coordinate subspace equals the intersection of $B$ with that subspace. Define the axial intercepts of $B$ by
$$b_i=\max\{t\ge0:t\mathbf{e}_i\in B\}=\max_{\mathbf{x}\in B}x_i>0.$$
After permuting coordinates, we assume $b_1\ge\cdots\ge b_n$.

For a convex body $K$, denote by $\operatorname{sh}_i K$ its coordinate Blaschke shaking~\cite{ref06} in direction $\mathbf{e}_i$. Over each point $\mathbf{y}$ of the projection onto $H_i=\mathbf{e}_i^\perp$, this operation replaces the fiber $[\alpha(\mathbf{y}),\beta(\mathbf{y})]$ in the $i$th coordinate by $[0,\beta(\mathbf{y})-\alpha(\mathbf{y})]$. Freyer and Lucas~\cite[Theorem 3.1]{ref09} construct the anti-blocking reduction by an initial unimodular transformation (i.e., an integer matrix with determinant $\pm1$) of coordinates followed by $\operatorname{sh}_1,\ldots,\operatorname{sh}_n$. This reduction preserves volume and does not decrease the lattice-point count.

\begin{theorem}[Freyer--Lucas]\label{thm:2-1}
For every convex body $K\subset\R^n$, there is an anti-blocking body $B\subset\R_{\ge0}^n$ such that
$$\vol_n(B)=\vol_n(K),\qquad \lambda_i(B_s)\le\lambda_i(K_s)\quad(i\in[n]),$$
$$\mathrm{G}(B)\ge\mathrm{G}(K),\qquad \mathrm{G}(\interior B)\le\mathrm{G}(\interior K).$$
\end{theorem}

For anti-blocking bodies, Freyer and Lucas \cite[Lemma 3.3]{ref09} also identify the successive minima with the ordered axial intercepts.

\begin{lemma}[Freyer--Lucas]\label{lem:2-2}
Let $B\subset\R_{\ge0}^n$ be an anti-blocking body, with axial intercepts $b_1\ge\cdots\ge b_n>0$. Then for all $i\in[n]$, one has
$$\lambda_i(B_s)=\frac{2}{b_i}.$$
\end{lemma}

\section{Compression toward earlier coordinates}\label{sec:3}

Inspired by the planar compression in \cite[Lemma 5.1]{ref09}, we introduce diagonal compressions in the coordinate planes. Maximizing a common linear moment over a compact admissible class gives a body fixed by every such compression, hence fully shifted, with the same volume and the required lattice-point inequalities. Fully shifted means that transferring part of a later coordinate to an earlier one leaves the point in the body. This allows us to place lattice cells below and above lattice points in the cumulative coordinates of Sect.~\ref{sec:4} and Sect.~\ref{sec:7}.

\subsection{A single diagonal compression}
Let $B\subset\R_{\ge0}^n$ be an anti-blocking body and fix $i<j$. Write
$$s=x_i+x_j,\qquad t=x_j,\qquad\mathbf{z}=(x_h)_{h\notin\{i,j\}}.$$
This is a unimodular transformation that preserves $\Z^n$. Let $\Pi$ be the projection of $B$ onto the $(s,\mathbf{z})$ coordinates. For $\mathbf{q}=(s,\mathbf{z})\in\Pi$, the corresponding closed fiber is an interval $[\alpha(\mathbf{q}),\beta(\mathbf{q})]$ with length $\ell(\mathbf{q})=\beta(\mathbf{q})-\alpha(\mathbf{q})$. Define $C_{ij}B$ by replacing this interval by $[0,\ell(\mathbf{q})]$, for every $\mathbf{q}\in\Pi$, and then returning to the original coordinates. In particular, a one-point fiber is retained as a one-point fiber. Geometrically, the operation moves each chord parallel to $\mathbf{e}_j-\mathbf{e}_i$ until its lower endpoint lies in $x_j=0$.

\begin{lemma}\label{lem:3-1}
The set $C_{ij}B$ is an anti-blocking body, and
$$\vol_n(C_{ij}B)=\vol_n(B),\qquad\mathrm{G}(C_{ij}B)\ge\mathrm{G}(B),$$
$$\mathrm{G}(\interior(C_{ij}B))\le\mathrm{G}(\interior B).$$
If $a_h(L)=\max\{r\ge0:r\mathbf{e}_h\in L\}$ denotes an axial intercept, then
\begin{equation}\label{eq:3-1}
 \begin{aligned}
 a_i(C_{ij}B)&=\max_{\mathbf{x}\in B}(x_i+x_j),\\
 a_j(C_{ij}B)&=\min\{a_i(B),a_j(B)\},
 \end{aligned}
\end{equation}
and the other axial intercepts are unchanged. The operation also preserves set inclusion.
\end{lemma}

\begin{proof}
In the $(s,t,\mathbf{z})$ coordinates, the operation of replacing $[\alpha(\mathbf{q}),\beta(\mathbf{q})]$ by $[0,\ell(\mathbf{q})]$ for every $\mathbf{q}\in\Pi$ is in fact a coordinate Blaschke shaking. Compactness, convexity and inclusion monotonicity follow from \cite[Lemma 1.1]{ref07}. Volume preservation and the lattice-point comparison follow from \cite[Lemma 3.2(ii), (v), (vi)]{ref09}, since the coordinate transformation preserves both volume and $\Z^n$.

It remains to prove that $C_{ij}B$ is an anti-blocking body and satisfies~\eqref{eq:3-1}. Fix $\mathbf{z}$ such that the section $B_{\mathbf{z}}=B\cap\{\mathbf{x}:x_h=z_h\text{ for }h\notin\{i,j\}\}$ is nonempty. We identify this section with a subset of $\R^2$ through the coordinates $(x_i,x_j)$. It is an anti-blocking convex set, possibly of lower dimension. Let $u,v$ be its axial intercepts and put $H=\max\{r+t:(r,t)\in B_{\mathbf{z}}\}$. For $0\le s\le H$, the second coordinate $x_j$ of $B_{\mathbf{z}}\cap\{\mathbf{x}:x_i+x_j=s\}$ ranges over $[\alpha,\beta]\subseteq[0,s]$. Its $x_j$-coordinate interval has length $\ell(s)=\beta-\alpha$. For $0\le s'\le s$, decreasing the coordinates of the two endpoints gives points at level $s'$ with second coordinates $\max(0,\alpha-s+s')$ and $\min(s',\beta)$, respectively. By convexity, we have
$$\ell(s')\ge\min(s',\beta)-\max(0,\alpha-s+s')\ge\min(s',\beta-\alpha).$$
Thus
\begin{equation}\label{eq:3-2}
 \ell(s')\ge\min\{s',\ell(s)\}\qquad(0\le s'\le s\le H).
\end{equation}
Let $(x,y)$ belong to the compressed section. Then $y\le\ell(x+y)\le\min\{u,v\}$. For $0\le x'\le x$ and $0\le y'\le y$, by~\eqref{eq:3-2}, we have
$$y'\le\min\{x'+y',\ell(x+y)\}\le\ell(x'+y').$$
Thus $(x',y')$ belongs to the compressed section. Its first intercept is $H$, since the compressed section contains the entire projection interval on the earlier axis. Its second intercept is $\min(u,v)$. Indeed, $\ell(s)=s$ holds exactly when the original fiber is $[0,s]$, or equivalently when both axial points $(s,0)$ and $(0,s)$ belong to the section. This argument also covers zero intercepts and one-point fibers.

For the remaining coordinates, the original sections satisfy $B_{\mathbf{z}}\subseteq B_{\mathbf{z}'}$ whenever $\mathbf{0}\le\mathbf{z}'\le\mathbf{z}$. Since compression preserves inclusion, we have $(C_{ij}B)_{\mathbf{z}}\subseteq(C_{ij}B)_{\mathbf{z}'}$. Let $\mathbf{x}\in C_{ij}B$ and $\mathbf{0}\le\mathbf{y}\le\mathbf{x}$, and write $\mathbf{z}=(x_h)_{h\notin\{i,j\}}$ and $\mathbf{z}'=(y_h)_{h\notin\{i,j\}}$. Inclusion monotonicity gives $(x_i,x_j)\in(C_{ij}B)_{\mathbf{z}'}$. Since this section is anti-blocking, it also contains $(y_i,y_j)$. Hence $\mathbf{y}\in C_{ij}B$, so $C_{ij}B$ is an anti-blocking body. The two selected intercepts are read from the section $\mathbf{z}=\mathbf{0}$. Setting all other coordinates to zero shows that~\eqref{eq:3-1} holds. The projection onto the remaining coordinates is unchanged, and every nonempty compressed planar section contains its origin. Hence the other intercepts are unchanged.
\end{proof}

\begin{figure}[!htbp]
\centering
\begin{tikzpicture}[scale=.64]
 \begin{scope}
  \fill[gray!12] (0,0) rectangle (3,2);
  \draw[thick] (0,0) rectangle (3,2);
  \draw[->] (-.3,0)--(5.5,0) node[right] {$x_i$};
  \draw[->] (0,-.3)--(0,2.7) node[above] {$x_j$};
  \draw[very thick] (2,2)--(3,1);
  \draw[densely dashed] (3,1)--(4,0);
  \fill (2,2) circle (2pt); \fill (3,1) circle (2pt);
  \node[below] at (3,0) {$3$}; \node[left] at (0,2) {$2$};
  \node at (2.5,-1.25) {original rectangle};
 \end{scope}
 \begin{scope}[xshift=8cm]
  \fill[gray!12] (0,0)--(5,0)--(1,2)--(0,2)--cycle;
  \draw[thick] (0,0)--(5,0)--(1,2)--(0,2)--cycle;
  \draw[->] (-.3,0)--(5.5,0) node[right] {$x_i$};
  \draw[->] (0,-.3)--(0,2.7) node[above] {$x_j$};
  \draw[densely dashed] (2,2)--(3,1);
  \draw[very thick] (3,1)--(4,0);
  \fill (3,1) circle (2pt); \fill (4,0) circle (2pt);
  \node[below] at (5,0) {$5$}; \node[left] at (0,2) {$2$};
  \node at (2.5,-1.25) {compressed trapezoid};
 \end{scope}
\end{tikzpicture}
\caption{Diagonal compression preserves each chord length. On the line
$x_i+x_j=4$, the solid segment moves from $1\le x_j\le2$ to
$0\le x_j\le1$. The same operation is performed on every parallel chord.}
\label{fig:1}
\end{figure}
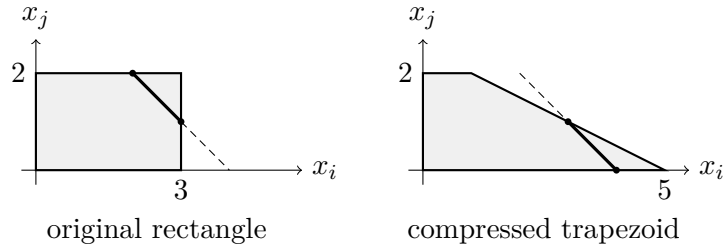

For example, $C_{12}$ sends $[0,3]\times[0,2]$ to the trapezoid with vertices $(0,0),(5,0),(1,2),(0,2)$; see Figure~\ref{fig:1}. The area remains $6$, while the larger axial intercept increases from $3$ to $5$. In general, a selected intercept pair $\{u,v\}$, with $u\le v$, is replaced by $\{u,H\}$, where $H\ge v$. It follows that the decreasingly ordered list of all intercepts can only increase componentwise.
\FloatBarrier

\subsection{A common fixed body}

\begin{definition}\label{def:3-2}
An anti-blocking body $A\subset\R_{\ge0}^n$ is {\it fully shifted} if
$$\mathbf{x}+\tau(\mathbf{e}_i-\mathbf{e}_j)\in A \quad\text{for }\mathbf{x}\in A,\ i<j,\ 0\le\tau\le x_j.$$
\end{definition}
This means that any amount of a later coordinate may be transferred to an earlier one. It is equivalent to $C_{ij}A=A$ for every $i<j$, because each nonempty diagonal fiber must then start at $x_j=0$. In particular, the axial intercepts of $A$ satisfy $a_1\ge\cdots\ge a_n$. The simplex $\Delta_n$ is a fully shifted body, whereas a positive coordinate box is not a fully shifted body when $n\ge2$.

\begin{proposition}\label{prop:3-3}
For every anti-blocking body $B\subset\R_{\ge0}^n$, there is a fully shifted body $A$ of the same volume such that
$$\mathrm{G}(A)\ge\mathrm{G}(B),\qquad \mathrm{G}(\interior A)\le\mathrm{G}(\interior B).$$
Its axial intercepts satisfy $a_i(A)\ge b_i$ for $1\le i\le n$, where $b_1\ge\cdots\ge b_n$ are the ordered intercepts of $B$.
\end{proposition}

\begin{proof}
For $n=1$, take $A=B$. Let $n\ge2$ and set $V=\vol_n(B),\eta=b_n>0$. For an anti-blocking body $L$, write $\hat a_1(L),\ldots,\hat a_n(L)$ for its intercepts arranged in decreasing order. Consider the nonempty class
$$\mathcal{F}=\left\{L\subset\R_{\ge0}^n: \begin{array}{l} L\text{ is an anti-blocking body},\ \vol_n(L)=V,\\ \mathrm{G}(L)\ge\mathrm{G}(B),\quad\mathrm{G}(\interior L)\le\mathrm{G}(\interior B),\\ \hat a_i(L)\ge b_i\quad(1\le i\le n) \end{array}\right\}.$$
We first show that $\mathcal{F}$ is compact. Every $L\in\mathcal{F}$ contains its axial simplex and has all intercepts at least $\eta$. Hence
$$\prod_{i=1}^n a_i(L)\le n!V,\qquad \eta\Delta_n\subseteq L\subseteq[0,J]^n,\qquad J=\frac{n!V}{\eta^{n-1}}.$$
By Blaschke's selection theorem, every sequence in $\mathcal{F}$ has a Hausdorff-convergent subsequence, again denoted by $(L_m)$, with limit $L^*$. The limit is anti-blocking and contains $\eta\Delta_n$. Volume is continuous under this convergence. The ordered axial intercepts are continuous as well, since they are rearrangements of the coordinate support values. Moreover, points outside $L^*$ eventually lie outside $L_m$, whereas points in $\interior L^*$ eventually lie in $\interior L_m$. The common box contains only finitely many lattice points, so
$$\begin{aligned}
\mathrm{G}(L^*)&\ge\limsup_m\mathrm{G}(L_m)\ge\mathrm{G}(B),\\
\mathrm{G}(\interior L^*)&\le\liminf_m\mathrm{G}(\interior L_m)\le\mathrm{G}(\interior B).
\end{aligned}$$
Thus $L^*\in\mathcal{F}$, proving compactness.

By Lemma~\ref{lem:3-1} and the intercept comparison, we have $C_{ij}\mathcal{F}\subseteq\mathcal{F}$ for every $i<j$. The continuous functional
$$\Phi(L)=\int_L\sum_{h=1}^n(n+1-h)x_h\,d\mathbf{x}$$
attains a maximum at some $A\in\mathcal{F}$.

For $L\in\mathcal{F}$, let $\Pi$, $\alpha$, and $\ell$ denote its projection and fiber data from the definition of $C_{ij}$. In the coordinates of a diagonal compression, the integrand is $(n+1-i)s-(j-i)t$ plus terms independent of $t$. Therefore
\begin{equation}\label{eq:3-4}
 \Phi(C_{ij}L)-\Phi(L)=(j-i)\int_\Pi
 \alpha(\mathbf{q})\ell(\mathbf{q})\,d\mathbf{q}\ge0.
\end{equation}
If equality holds, then $\alpha=0$ on almost every positive-length fiber. Thus $C_{ij}L$ and $L$ agree almost everywhere and, being full-dimensional compact convex bodies, are equal.

Since $C_{ij}A\in\mathcal{F}$, maximality and~\eqref{eq:3-4} force $C_{ij}A=A$ for every $i<j$. Thus, $A$ is a fully shifted body. Its volume, two lattice-point comparisons and ordered intercept bounds follow from $A\in\mathcal{F}$; the intercepts are now already in coordinate order.
\end{proof}

Together with Theorem~\ref{thm:2-1} and Lemma~\ref{lem:2-2}, Proposition~\ref{prop:3-3} reduces the upper bound~\eqref{eq:1-3} in Conjecture~\ref{conj:1-1} to
$$\mathrm{G}(A)\le\vol_n(A)\prod_{i=1}^n\left(1+\frac{i}{a_i}\right)$$
for a fully shifted body with intercepts $a_1\ge\cdots\ge a_n>0$. For the lower bound~\eqref{eq:1-4} in Conjecture~\ref{conj:1-1}, the same reduction gives the target
$$\mathrm{G}(\interior A)\ge\vol_n(A) \prod_{i=1}^n\left(1-\frac{i}{a_i}\right),\quad\text{if }a_n\ge n.$$
\section{Lattice cells and section volumes}\label{sec:4}

Let $A\subset\R_{\ge0}^n$ be a fully shifted body, as in Definition~\ref{def:3-2}. In each coordinate section, we pass to cumulative coordinates and place cells below lattice points. A generating-function identity will combine their estimates into a bound for $\mathrm{G}(A)$.

Let $a_1\ge\cdots\ge a_n>0$ be the axial intercepts of $A$. For each nonempty set $S\subseteq[n]$, put
$$E_S=\operatorname{span}\{\mathbf{e}_i:i\in S\},\qquad A_S=A\cap E_S,\qquad V_S=\vol_{|S|}(A_S),$$
where volume is taken in the coordinate subspace $E_S$. The section $A_S$ contains its axial simplex, so it is full-dimensional in $E_S$. It inherits the anti-blocking property and every allowed transfer between coordinates in $S$. For $T\subseteq[n]$, define
$$N_T=\#\{\mathbf{x}\in A\cap\Z^n:\{i:x_i>0\}=T\}.$$
Thus $N_\varnothing=1$ and $\mathrm{G}(A)=1+\sum_{\varnothing\ne T\subseteq[n]}N_T$.

\begin{lemma}\label{lem:4-1}
Let $S=\{s_1<\cdots<s_d\}\subseteq[n]$ be a nonempty set. For a nonempty set $T\subseteq S$ with $\max T=\max S$, write $T=\{s_{r_1},\ldots,s_{r_\ell}\}$, where $r_1<\cdots<r_\ell=d$, and set
$$c(S,T)=\prod_{h=1}^\ell\frac{1}{(r_h-r_{h-1})!},\qquad r_0=0.$$
Set $c(S,T)=0$ for all other nonempty sets $T\subseteq[n]$. Then
\begin{equation}\label{eq:4-2}
V_S\ge\sum_{\varnothing\ne T\subseteq[n]}c(S,T)N_T.
\end{equation}
\end{lemma}

\begin{proof}
We use the cumulative coordinates on $E_S$
$$y_j=x_{s_j}+\cdots+x_{s_d}\qquad(1\le j\le d).$$
This is a unimodular transformation from $E_S\cap\Z^n$ to $\Z^d$, and it preserves $d$-dimensional volume. Let $D_S$ be the image of $A_S$. It lies in the ordered cone
$$W_d=\{\mathbf{y}\in\R^d:y_1\ge\cdots\ge y_d\ge0\}.$$
We claim that $D_S$ is an order ideal with respect to $W_d$, meaning that $\mathbf{q}\in D_S$ whenever $\mathbf{y}\in D_S$, $\mathbf{q}\in W_d$, and $\mathbf{q}\le\mathbf{y}$. For $d=1$, this is the anti-blocking property. Suppose $d\ge2$ and put $q_{d+1}=0$. Transfer $y_j-q_j$ from $x_{s_j}$ to $x_{s_{j-1}}$ in the order $j=d,d-1,\ldots,2$. Finally, decrease the first coordinate by $y_1-q_1$, using the anti-blocking property. The resulting coordinates are
\begin{gather*}
x^*_{s_d}=y_d-q_{d+1}-(y_d-q_d)=q_d;\\ x^*_{s_{d-1}}=x_{s_{d-1}}+(y_d-q_d)-(y_{d-1}-q_{d-1})=q_{d-1}-q_d;\\ \cdots;\\ x^*_{s_j}=x_{s_j}+(y_{j+1}-q_{j+1})-(y_j-q_j)=q_j-q_{j+1}\qquad(2\le j\le d-1);\\ \cdots;\\ x^*_{s_1}=x_{s_1}+(y_2-q_2)-(y_1-q_1)=q_1-q_2.
\end{gather*}
Each transfer is permitted because $A_S$ is fully shifted and $0\le y_j-q_j\le y_j-q_{j+1}$. The final decrease leaves the first coordinate equal to $q_1-q_2\ge0$. Thus the resulting point $\mathbf{x}^*\in A_S$ satisfies
$$q_j=x^*_{s_j}+\cdots+x^*_{s_d}\qquad(1\le j\le d).$$
So we have $\mathbf{q}\in D_S$.

We estimate $V_S$ by constructing cells below lattice points in $D_S$. For each $\mathbf{y}\in D_S\cap\Z^d$ with $y_d\ge1$, we define the closed cell
$$Q^-_{\mathbf{y}}=W_d\cap\prod_{j=1}^d[y_j-1,y_j].$$
Because $D_S$ is an order ideal with respect to $W_d$, we have $Q^-_{\mathbf{y}}\subseteq D_S$. Distinct cells have disjoint interiors.

For a nonzero lattice point of $A_S$ with positive-coordinate support set $T=\{s_{r_1},\ldots,s_{r_\ell}\}$, the condition $y_d\ge1$ is equivalent to $\max T=\max S$. Moreover, if $y_j>y_{j+1}$, integer spacing gives $y_j-1\ge y_{j+1}$, so this order constraint holds throughout the box. The equal-coordinate blocks then have lengths $r_1,r_2-r_1,\ldots,r_\ell-r_{\ell-1}$. Since a block of $b$ equal coordinates occupies $1/b!$ of its unit box (due to the $b!$ symmetric orderings), the volume of the cell is exactly $c(S,T)$. Adding their volumes proves \eqref{eq:4-2}, including the contributions from lattice points on the boundary of $A$.
\end{proof}

For example, in dimension two, the coordinates are $y_1=x_1+x_2$ and $y_2=x_2$. For $y_1>y_2\ge1$ the cell is a full unit square; for $y_1=y_2\ge1$ it is half of that square. Thus
$$V_{\{1,2\}}\ge N_{\{1,2\}}+\tfrac{1}{2}N_{\{2\}}.$$
Points on the first axis have $y_2=0$ and contribute no area. The lower-dimensional sections must therefore supply the remaining contributions.

We now choose weights that give every nonzero lattice point total coefficient one. For $S=\{s_1<\cdots<s_d\}\ne\varnothing$, define
\begin{equation}\label{eq:4-3}
w_S=\prod_{j=1}^d\frac{1}{s_j-s_{j-1}},\qquad s_0=0,\quad w_\varnothing=1.
\end{equation}
The gaps in \eqref{eq:4-3} are gaps between the original coordinate indices. The factorials in the definition of $c(S,T)$ use positions within $S$. In particular, $w_{[n]}=1$ and $w_{\{i\}}=1/i$. For example, when $n\ge4$, take $S=\{1,3,4\}$ and $T=\{4\}$. Then $w_S=1/(1\cdot2\cdot1)=1/2$, whereas $c(S,T)=1/3!=1/6$. Here the original indices have gaps $1,2,1$, but a point supported only on coordinate $4$ has three equal cumulative coordinates in this section.

\begin{proposition}\label{prop:4-2}
Every fully shifted anti-blocking body $A\subset\R_{\ge0}^n$ satisfies
\begin{equation}\label{eq:4-4}
\mathrm{G}(A)\le1+\sum_{\varnothing\ne S\subseteq[n]}w_SV_S.
\end{equation}
\end{proposition}

\begin{proof}
We first show that, for every nonempty set $T\subseteq[n]$,
\begin{equation}\label{eq:4-5}
\sum_{\varnothing\ne S\subseteq[n]}w_Sc(S,T)=1.
\end{equation}
Fix $T=\{t_1<\cdots<t_q\}$ and put $t_0=0$. Then $c(S,T)\ne0$ only if $S\supseteq T$ and $\max S=t_q$. Such a set $S$ is determined independently on the integer intervals $(t_{h-1},t_h]$, each of which must contain its right endpoint among the selected indices.

Consider an interval of length $b=t_h-t_{h-1}$. If $S$ has $\ell$ indices in this interval, their successive gaps form a composition
$$p_1+\cdots+p_\ell=b,\qquad p_j\in\Z_{>0}.$$
Its contribution to $w_S$ is $1/(p_1\cdots p_\ell)$, and its contribution to $c(S,T)$ is $1/\ell!$. Therefore, the total contribution of the interval is
$$\gamma_b=\sum_{\ell=1}^b\frac{1}{\ell!} \sum_{\substack{p_1+\cdots+p_\ell=b\\p_j\ge1}}\frac{1}{p_1\cdots p_\ell}.$$
The generating function of these coefficients is
$$1+\sum_{b\ge1}\gamma_bz^b =\exp\left(\sum_{p\ge1}\frac{z^p}{p}\right)=\frac{1}{1-z}.$$
Hence $\gamma_b=1$ for every positive integer $b$. The choices on the intervals are independent, so the left side of \eqref{eq:4-5} is $\prod_{h=1}^q\gamma_{t_h-t_{h-1}}=1$.

Multiplying \eqref{eq:4-2} by $w_S$, summing over all nonempty $S$, and using \eqref{eq:4-5}, we obtain
\begin{align*}
\mathrm{G}(A)&=1+\sum_{\varnothing\ne T\subseteq[n]}N_T\\
&=1+\sum_{\varnothing\ne T\subseteq[n]}N_T
\sum_{\varnothing\ne S\subseteq[n]}w_Sc(S,T)\\
&=1+\sum_{\varnothing\ne S\subseteq[n]}w_S
\sum_{\varnothing\ne T\subseteq[n]}c(S,T)N_T\\
&\le1+\sum_{\varnothing\ne S\subseteq[n]}w_SV_S.
\end{align*}
\end{proof}

On a one-dimensional section, Lemma~\ref{lem:4-1} reads $a_i\ge\lfloor a_i\rfloor$. Thus the axial terms in \eqref{eq:4-4} are the real lengths $a_i/i$; the estimate applies to arbitrary positive axial intercepts. To compare section volumes with the full volume, we successively join a section to the complementary axial endpoints.

\begin{lemma}\label{lem:4-3}
Let $A\subset\R_{\ge0}^n$ be an anti-blocking body with axial intercepts $a_1,\ldots,a_n$, and set $A_\varnothing=\{\mathbf{0}\}$ and $V_\varnothing=1$. Then, for every set $S\subseteq[n]$,
$$V_S\le\frac{n!}{|S|!}\vol_n(A)\prod_{j\notin S}a_j^{-1}.$$
\end{lemma}

\begin{proof}
Let $k=|S|$ and set
$$R=\conv\bigl(A_S\cup\{a_j\mathbf{e}_j:j\notin S\}\bigr)\subseteq A.$$
Joining $A_S$ successively to the complementary axial endpoints gives pyramids of heights $a_j$ and dimensions $k+1,\ldots,n$. Repeated use of the pyramid volume formula, with $\vol_0(\{\mathbf{0}\})=1$, gives
$$\vol_n(R)=\frac{k!}{n!}V_S\prod_{j\notin S}a_j.$$
The inclusion $R\subseteq A$ proves the assertion.
\end{proof}

For a fully shifted body $A$, Proposition~\ref{prop:4-2} and Lemma~\ref{lem:4-3} therefore give
\begin{equation}\label{eq:4-6}
\begin{aligned}
\mathrm{G}(A)&\le\sum_{S\subseteq[n]}w_SV_S\\
&\le\frac{n!\vol_n(A)}{\prod_{i=1}^na_i}
\sum_{S\subseteq[n]}\frac{w_S}{|S|!}\prod_{i\in S}a_i.
\end{aligned}
\end{equation}
It remains to compare the polynomial on the right with $\prod_{i=1}^n(1+a_i/i)$, using the order $a_1\ge\cdots\ge a_n$.

\section{The ordered-sequence inequality}\label{sec:5}

We bound the subset sum in \eqref{eq:4-6} by $\prod_{i=1}^n(1+a_i/i)$. An elementary induction gives a comparison in each degree.

\begin{proposition}\label{prop:5-1}
Let $n\ge1$ and $a_1\ge\cdots\ge a_n\ge0$. For every $0\le k\le n$,
\begin{equation}\label{eq:5-1}
\frac{1}{k!}\sum_{1\le s_1<\cdots<s_k\le n}
\prod_{r=1}^k\frac{a_{s_r}}{s_r-s_{r-1}}
\le\sum_{1\le s_1<\cdots<s_k\le n}
\prod_{r=1}^k\frac{a_{s_r}}{s_r},\qquad s_0=0.
\end{equation}
For $k=0$, both sides are understood to be $1$. For $k\ge1$, the inequality also holds when both sums are restricted to $s_k=m$, for any $k\le m\le n$.
\end{proposition}

\begin{proof}
The case $k=0$ follows from the stated convention. Fix $1\le k\le m\le n$ and let
$$\mathcal{C}_{m,k}=\{(g_1,\ldots,g_k)\in\Z_{>0}^k:g_1+\cdots+g_k=m\}.$$
Taking partial sums gives a bijection between $\mathcal{C}_{m,k}$ and the index sets $1\le s_1<\cdots<s_k=m$. For $\mathbf{g}\in\mathcal{C}_{m,k}$ and a label permutation $\pi\in\mathfrak{S}_k$, put
$$p_r(\pi)=g_{\pi_1}+\cdots+g_{\pi_r},\qquad F_{\mathbf{g}}(\pi)=\prod_{r=1}^ka_{p_r(\pi)}.$$
All indices $p_r(\pi)$ lie between $1$ and $m$, and equal parts keep distinct labels. Define
$$\mu(\mathbf{g})=\frac{1}{k!}\sum_{\pi\in\mathfrak{S}_k}F_{\mathbf{g}}(\pi),\qquad \nu(\mathbf{g})=\left(\prod_{i=1}^kg_i\right) \sum_{\pi\in\mathfrak{S}_k}\frac{F_{\mathbf{g}}(\pi)}{\prod_{r=1}^kp_r(\pi)}.$$
We prove $\nu(\mathbf{g})\ge \mu(\mathbf{g})$ for all positive integer tuples of total at most $n$, by induction on their number of parts. For $k=1$, both expressions equal $a_{g_1}$. For $k\ge2$, write $\mu_j=\mu(\mathbf{g}^{\widehat{j}})$ and $\nu_j=\nu(\mathbf{g}^{\widehat{j}})$, where $\mathbf{g}^{\widehat{j}}$ is obtained by deleting label $j$. Grouping permutations by their last label and using $p_k(\pi)=m$ gives
$$\mu(\mathbf{g})=\frac{a_m}{k}\sum_{j=1}^k\mu_j,\qquad \nu(\mathbf{g})=\frac{a_m}{m}\sum_{j=1}^kg_j\nu_j.$$
If $g_j\ge g_\ell$, replace label $\ell$ by label $j$ in each permutation of the labels other than $j$. This gives a bijection to permutations of the labels other than $\ell$. Each prefix sum is unchanged before the replaced position and increases by $g_j-g_\ell$ thereafter. Since the $a_i$ are nonincreasing and nonnegative, the corresponding product cannot increase. Hence $\mu_j\ge \mu_\ell$, and therefore
$$\sum_{j=1}^kg_j\mu_j-\frac{m}{k}\sum_{j=1}^k\mu_j =\frac{1}{k}\sum_{j<\ell}(g_j-g_\ell)(\mu_j-\mu_\ell)\ge0.$$
The induction hypothesis now yields
$$\nu(\mathbf{g})\ge\frac{a_m}{m}\sum_jg_j\mu_j \ge\frac{a_m}{k}\sum_j\mu_j=\mu(\mathbf{g}).$$
Dividing this comparison by $\prod_ig_i$ gives
\begin{equation}\label{eq:5-2}
\sum_{\pi\in\mathfrak{S}_k}\frac{F_{\mathbf{g}}(\pi)}{\prod_{r=1}^kp_r(\pi)}
\ge\frac{1}{k!}\sum_{\pi\in\mathfrak{S}_k}
\frac{F_{\mathbf{g}}(\pi)}{\prod_{i=1}^kg_i}.
\end{equation}
Sum \eqref{eq:5-2} over all $\mathbf{g}\in\mathcal{C}_{m,k}$. For each fixed label permutation $\pi$, the map $\mathbf{g}\mapsto(g_{\pi_1},\ldots,g_{\pi_k})$ is a bijection of $\mathcal{C}_{m,k}$, including compositions with repeated parts. Consequently the sum becomes
$$k!\sum_{\mathbf{g}\in\mathcal{C}_{m,k}} \frac{\prod_{r=1}^ka_{g_1+\cdots+g_r}}{\prod_{r=1}^k(g_1+\cdots+g_r)} \ge\sum_{\mathbf{g}\in\mathcal{C}_{m,k}}\prod_{r=1}^k\frac{a_{g_1+\cdots+g_r}}{g_r}.$$
Divide by $k!$ and return to the partial sums $s_r=g_1+\cdots+g_r$. This proves \eqref{eq:5-1} with $s_k=m$. Summing over $m=k,\ldots,n$ proves the proposition.
\end{proof}

\begin{corollary}\label{cor:5-2}
For $n\ge1$ and $a_1\ge\cdots\ge a_n\ge0$, the weights in \eqref{eq:4-3} satisfy
\begin{equation}\label{eq:5-3}
\sum_{S\subseteq[n]}\frac{w_S}{|S|!}\prod_{i\in S}a_i
\le\prod_{i=1}^n\left(1+\frac{a_i}{i}\right).
\end{equation}
\end{corollary}

\begin{proof}
Summing \eqref{eq:5-1} over $k=0,\ldots,n$, we obtain
$$
\begin{aligned}
\sum_{S\subseteq[n]}\frac{w_S}{|S|!}\prod_{i\in S}a_i
&=\sum_{k=0}^n\sum_{\substack{S\subseteq[n]\\|S|=k}}\frac{w_S}{k!}\prod_{i\in S}a_i\\
&\le\sum_{S\subseteq[n]}\prod_{i\in S}\frac{a_i}{i}
=\prod_{i=1}^n\left(1+\frac{a_i}{i}\right).
\end{aligned}
$$
\end{proof}

\section{Proof of Theorem~\ref{thm:1-1}}\label{sec:6}

\subsection{The inequality part}
\begin{proof}
By Theorem~\ref{thm:2-1}, Lemma~\ref{lem:2-2} and Proposition~\ref{prop:3-3}, there is an anti-blocking body $B$ and a fully shifted body $A$ such that
$$
\begin{gathered}
\vol_n(A)=\vol_n(B)=\vol_n(K),\qquad \mathrm{G}(A)\ge \mathrm{G}(B)\ge \mathrm{G}(K),\\
a_i\ge b_i,\qquad\frac{2}{b_i}=\lambda_i(B_s)\le\lambda_i(K_s),
\end{gathered}
$$
where $a_1\ge\cdots\ge a_n>0$ and $b_1\ge\cdots\ge b_n>0$ are the ordered intercepts of $A$ and $B$, respectively. By \eqref{eq:4-6} and Corollary~\ref{cor:5-2}, we have
\begin{align*}
\mathrm{G}(K)\le \mathrm{G}(B)\le \mathrm{G}(A)&\le\frac{n!\vol_n(A)}{\prod_i a_i}
\prod_{i=1}^n\left(1+\frac{a_i}{i}\right)\\
&=\vol_n(A)\prod_{i=1}^n\left(1+\frac{i}{a_i}\right)
\le\vol_n(B)\prod_{i=1}^n\left(1+\frac{i}{b_i}\right)\\
&=\vol_n(B)\prod_{i=1}^n\left(1+\frac{i\lambda_i(B_s)}{2}\right)\\
&\le\vol_n(K)\prod_{i=1}^n\left(1+\frac{i\lambda_i(K_s)}{2}\right).
\end{align*}
\end{proof}

\subsection{The equality part}
We first determine equality for anti-blocking bodies. We then recover the original body by reversing the coordinate Blaschke shakings in the Freyer--Lucas reduction. An affine lattice transformation is a map $\mathbf{x}\mapsto U\mathbf{x}+\mathbf{z}$ with $U\in\GL_n(\Z)$ and $\mathbf{z}\in\Z^n$. The lattice-point equality will force each reversed shaking to have this form.

\begin{lemma}\label{lem:6-1}
Let $B\subset\R_{\ge0}^n$ be an anti-blocking body with ordered axial intercepts $b_1\ge\cdots\ge b_n>0$. Then
$$\mathrm{G}(B)=\vol_n(B)\prod_{i=1}^n\left(1+\frac{i}{b_i}\right)$$
if and only if $B=m\Delta_n$ for a positive integer $m$.
\end{lemma}

\begin{proof}
Suppose equality holds. Proposition~\ref{prop:3-3} gives a fully shifted $A$ of the same volume $V$, with $\mathrm{G}(A)\ge \mathrm{G}(B)$ and ordered intercepts $a_i\ge b_i$. Hence
$$\mathrm{G}(B)\le \mathrm{G}(A)\le V\prod_i\left(1+\frac{i}{a_i}\right) \le V\prod_i\left(1+\frac{i}{b_i}\right)=\mathrm{G}(B).$$
Since every comparison is an equality, the strict monotonicity of the factors implies $a_i=b_i$ for every $i$. Equality also holds at every step of \eqref{eq:4-6}. Since all weights are positive, for every set $S\subseteq[n]$,
\begin{equation}\label{eq:6-1}
V_S=\frac{n!}{|S|!}V\prod_{j\notin S}a_j^{-1}.
\end{equation}
Take $S=\varnothing$ in \eqref{eq:6-1}, we get
$$V=\frac{\prod_i a_i}{n!}.$$
In addition, both $A$ and $B$ contain the common axial simplex, so
$$A=B=\conv(\mathbf{0},a_1\mathbf{e}_1,\ldots,a_n\mathbf{e}_n).$$
Next, we prove that all intercepts $a_i$ are equal integers. Using \eqref{eq:4-5}, we have
$$0=\sum_{S\subseteq[n]}w_SV_S-\mathrm{G}(A) =\sum_{\varnothing\ne S\subseteq[n]}w_S \left(V_S-\sum_{\varnothing\ne T\subseteq[n]}c(S,T)N_T\right).$$
Every summand is nonnegative by Lemma~\ref{lem:4-1} and every weight is positive. Hence, for each nonempty set $S\subseteq[n]$, we have
\begin{equation}\label{eq:6-2}
V_S=\sum_{\varnothing\ne T\subseteq[n]}c(S,T)N_T.
\end{equation}
If $n=1$, the singleton estimate gives $a_1=\lfloor a_1\rfloor$, as required. Let $n\ge2$ and consider $S=\{1,n\}$. Put $p=a_1$, $q=a_n$ and $h=\lfloor q\rfloor$, so $p\ge q>0$. The positive lattice points on the last axis give $N_{\{n\}}=h$. For each $1\le j\le h$, the positive first coordinates in the slice $x_n=j$ range from $1$ to $\lfloor p(1-j/q)\rfloor$. Thus
$$N_{\{1,n\}}=\sum_{j=1}^h\left\lfloor p\left(1-\frac{j}{q}\right)\right\rfloor.$$
Taking $S=\{1,n\}$ in \eqref{eq:6-2}, only $T=\{1,n\}$ and $T=\{n\}$ contribute. Hence
\begin{align*}
0&=V_{\{1,n\}}-\sum_{\varnothing\ne T\subseteq[n]}c(\{1,n\},T)N_T\\
&=V_{\{1,n\}}-c(\{1,n\},\{1,n\})N_{\{1,n\}}-c(\{1,n\},\{n\})N_{\{n\}}\\
&=\frac{pq}{2}-\sum_{j=1}^h\left\lfloor p\left(1-\frac{j}{q}\right)\right\rfloor-\frac{h}{2}\\
&\ge\frac{p(q-h)^2+(p-q)h}{2q}\ge0.
\end{align*}
Then $p=q=h$, so all intercepts equal the positive integer $m=q$ and $B=m\Delta_n$. Conversely, for $B=m\Delta_n$ with $m\in\Z_{>0}$, counting nonnegative integer coordinates gives
$$\mathrm{G}(B)=\binom{m+n}{n}=\vol_n(B)\prod_{i=1}^n\left(1+\frac{i}{m}\right).$$
\end{proof}

\begin{lemma}\label{lem:6-2}
Let $K\subset\R^n$ be a convex body and let $P=\operatorname{sh}_i K$. Suppose that $P$ is an affine lattice image of $m\Delta_n$, where $m\in\Z_{>0}$, and that $\mathrm{G}(K)=\mathrm{G}(P)$. Then $K$ is also an affine lattice image of $m\Delta_n$.
\end{lemma}

\begin{proof}
For $n=1$, an interval of length $m$ contains $m+1$ integers precisely when its endpoints are integers. Assume $n\ge2$. We first determine the fibers of $P$. Since $P=\operatorname{sh}_i K$, $P\subseteq\{\mathbf{x}:x_i\ge0\}$ and $P\cap H_i$ is a facet of the simplex $P$. Put $P_i=P\cap H_i$, label its vertices $\mathbf{v}_0,\ldots,\mathbf{v}_{n-1}$, and denote the remaining vertex of $P$ by $\mathbf{v}_n$. The vectors
$$\mathbf{u}_r=\frac{\mathbf{v}_r-\mathbf{v}_0}{m},\qquad1\le r\le n,$$
form a basis of $\Z^n$. For all $1\le r<n$, we have $\mathbf{u}_r\in H_i$. Expanding the determinant in the $i$th row gives
$$1=(\mathbf{u}_n)_i\left|\det_{H_i}(\mathbf{u}_1,\ldots,\mathbf{u}_{n-1})\right|.$$
Here the determinant in $H_i$ is taken in $H_i\cap\Z^n$. Both factors are positive integers, so $(\mathbf{u}_n)_i=1$ and $\mathbf{u}_1,\ldots,\mathbf{u}_{n-1}$ form a basis of $H_i\cap\Z^n$.

The projection $\mathbf{w}=\mathbf{v}_n-m\mathbf{e}_i$ belongs to the base facet $P_i$. Moreover,
$$\frac{\mathbf{w}-\mathbf{v}_0}{m}=\mathbf{u}_n-\mathbf{e}_i \in\conv(\mathbf{0},\mathbf{u}_1,\ldots,\mathbf{u}_{n-1})\cap\Z^n.$$
A unimodular simplex contains no lattice points other than its vertices, so $\mathbf{w}$ is one of the vertices of $P_i$. Relabel the base vertices so that $\mathbf{w}=\mathbf{v}_0$, and redefine the $\mathbf{u}_r$ accordingly. Then $\mathbf{u}_1,\ldots,\mathbf{u}_{n-1}$ still form a basis of $H_i\cap\Z^n$, and $\mathbf{u}_n=\mathbf{e}_i$. For $\mathbf{y}=\mathbf{v}_0+\sum_{r=1}^{n-1}t_r\mathbf{u}_r\in P_i$, the fiber length is
$$\ell(\mathbf{y})=m-\sum_{r=1}^{n-1}t_r.$$
Thus $\ell$ is affine on $P_i$ and is an integer at every $\mathbf{y}\in P_i\cap\Z^n$.

Write the fibers of $K$ over $P_i$ as $[\alpha(\mathbf{y}),\beta(\mathbf{y})]$. Convexity of $K$ makes $\alpha$ convex and $\beta$ concave on $P_i$. Since $\beta-\alpha=\ell$ is affine on $P_i$, the function $\alpha$ is both convex and concave on $P_i$, and hence affine.

For each $\mathbf{y}\in P_i\cap\Z^n$, the fiber of $P$ contains $\ell(\mathbf{y})+1$ integer points. An interval of the nonnegative integer length $\ell(\mathbf{y})$ has this number of integer points if and only if its left endpoint is an integer. The lattice-point count on each individual fiber does not decrease under shaking. Since $\mathrm{G}(K)=\mathrm{G}(P)$, the left endpoint $\alpha(\mathbf{y})$ of the fiber $[\alpha(\mathbf{y}),\beta(\mathbf{y})]$ over every $\mathbf{y}\in P_i\cap\Z^n$ is an integer.

Since $m\ge1$, the points $\mathbf{v}_0,\mathbf{v}_0+\mathbf{u}_1,\ldots,\mathbf{v}_0+\mathbf{u}_{n-1}$ lie in $P_i$. They form an affine basis of $H_i\cap\Z^n$. The affine function $\alpha$ is integer-valued on this basis, so its affine extension is integer-valued on all of $H_i\cap\Z^n$. Consequently
$$\mathbf{y}+t\mathbf{e}_i\longmapsto\mathbf{y}+(t+\alpha(\mathbf{y}))\mathbf{e}_i$$
is an affine lattice transformation of determinant one. It maps $P$ onto $K$, proving the assertion.
\end{proof}

\begin{proof}[Equality in Theorem~\ref{thm:1-1}]
Work first with $\Lambda=\Z^n$, and suppose equality holds in~\eqref{eq:1-3}. Every comparison in the proof of the upper bound is then an equality. In particular, the anti-blocking reduction $B$ attains its upper bound. Lemma~\ref{lem:6-1} gives $B=m\Delta_n$ for an integer $m\ge1$.

The coordinate permutation used to order the intercepts leaves $m\Delta_n$ unchanged. Thus the finite construction of the anti-blocking body gives, for some $U\in\GL_n(\Z)$,
$$K_0=UK,\qquad K_j=\operatorname{sh}_jK_{j-1}\quad(1\le j\le n),\qquad K_n=B.$$
Since 
$$\mathrm{G}(K)=\mathrm{G}(K_0)\le \mathrm{G}(K_1)\le\cdots\le \mathrm{G}(K_n)=\mathrm{G}(B)=\mathrm{G}(K),$$ 
we have $\mathrm{G}(K_{j-1})=\mathrm{G}(K_j)$ for all $1\le j\le n$. Starting with $K_n=m\Delta_n$, Lemma~\ref{lem:6-2} applies successively for $j=n,n-1,\ldots,1$. It follows that $K_0$ is an affine lattice image of $m\Delta_n$. By inverting the initial coordinate transformation and the lattice normalization, we recover the desired description for an arbitrary full-rank lattice.

Conversely, it follows directly from Lemma~\ref{lem:6-1} that every body of the stated form attains equality.
\end{proof}

\section{A lower bound for interior lattice points}\label{sec:7}

We use the cumulative coordinates from Sect.~\ref{sec:4} again, this time placing unit cells above rounded points. An upper cell with repeated coordinates at its anchor is assigned to a lower cell in a coordinate section of codimension one. The latter cell has at least twice its volume.

\begin{lemma}\label{lem:7-1}
Every fully shifted anti-blocking body $A\subset\R_{\ge0}^n$ satisfies
\begin{equation}\label{eq:7-1}
\mathrm{G}(\interior A)\ge\vol_n(A)-V_{[n]\setminus\{n\}}-\frac{1}{2}\sum_{i=1}^{n-1}V_{[n]\setminus\{i\}},
\end{equation}
where $V_\varnothing=1$.
\end{lemma}

\begin{proof}
For $n=1$, write $A=[0,a_1]$. Then $\mathrm{G}(\interior A)=\lceil a_1\rceil-1\ge a_1-1$. Suppose $n\ge2$, and let $D=D_{[n]}\subseteq W_n$ be the image of $A$ in the cumulative coordinates from Lemma~\ref{lem:4-1}. By that proof, $D$ is an order ideal in $W_n$.

Consider $\mathbf{y}\in D$ with $y_n>1$, excluding the measure-zero union of the hyperplanes $y_j\in\Z$. Define $\mathbf{q}=(q_1,\ldots,q_n)$ by $q_i=\lfloor y_i\rfloor$ for $1\le i\le n$. Then $\mathbf{q}\in D\cap\Z^n$, $q_n\ge1$, and $\mathbf{q}<\mathbf{y}$. Every $\mathbf{y}$ lies in the upper cell
$$Q^+_{\mathbf{q}}=W_n\cap\prod_{j=1}^n[q_j,q_j+1].$$
The block-volume of these upper cells can be calculated in the same way as Lemma~\ref{lem:4-1}.

If $q_1>\cdots>q_n>0$, then $\mathbf{q}$ has an open neighborhood in $\R^n$ contained in $W_n\cap\{\mathbf{z}:\mathbf{z}\le\mathbf{y}\}$. By the order-ideal property, this neighborhood lies in $D$, so $\mathbf{q}\in\interior D$. Thus $\mathbf{q}$ corresponds to an interior lattice point of $A$, and $Q^+_{\mathbf{q}}$ has volume one.

For each remaining anchor, choose the first repeated position
$$i=\min\{j<n:q_j=q_{j+1}\},\qquad S_i=[n]\setminus\{i\}.$$
Delete $q_i$ and denote the resulting vector by $\widehat{\mathbf{q}}$. The original-coordinate point has $x_i=q_i-q_{i+1}=0$, so it belongs to $A_{S_i}$. Its cumulative coordinates in that section are precisely $\widehat{\mathbf{q}}$, whose last coordinate is still $q_n\ge1$. Consequently, the cell $Q^-_{\widehat{\mathbf{q}}}$ from Lemma~\ref{lem:4-1} lies in the cumulative-coordinate image of $A_{S_i}$.

If the block containing $q_i$ has length $b\ge2$, deletion shortens that block to $b-1$ and leaves all other blocks unchanged. Therefore
$$\vol_{n-1}(Q^-_{\widehat{\mathbf{q}}})=b\vol_n(Q^+_{\mathbf{q}})\ge2\vol_n(Q^+_{\mathbf{q}}).$$
For fixed $i$, the deletion map is injective, since the removed entry must equal the following entry. Thus the lower cells assigned to $A_{S_i}$ are distinct and have disjoint interiors. Their total volume is at most $V_{S_i}$. Summing first within each section and then over $i<n$, and including the cells with strictly decreasing anchors, gives
\begin{equation}\label{eq:7-2}
\vol_n(D\cap\{y_n>1\})\le\mathrm{G}(\interior A)+\frac{1}{2}\sum_{i=1}^{n-1}V_{[n]\setminus\{i\}}.
\end{equation}

Since $y_n=x_n$, the remaining slab corresponds to $A\cap\{x_n\le1\}$, whose volume is at most $V_{[n]\setminus\{n\}}$ by the anti-blocking property and Fubini's theorem. Together with \eqref{eq:7-2}, this proves \eqref{eq:7-1}.
\end{proof}

\begin{proof}[Proof of Theorem~\ref{thm:1-2}]
Work with $\Z^n$ and put $V=\vol_n(K)$. By Theorem~\ref{thm:2-1}, Lemma~\ref{lem:2-2} and Proposition~\ref{prop:3-3}, there is a fully shifted body $A$ of volume $V$, with axial intercepts $a_1\ge\cdots\ge a_n>0$, such that
$$\mathrm{G}(\interior A)\le\mathrm{G}(\interior K),\qquad \frac{1}{a_i}\le\frac{\lambda_i(K_s)}{2}\quad(1\le i\le n).$$
Lemmas~\ref{lem:7-1} and~\ref{lem:4-3} therefore give
$$
\begin{aligned}
\mathrm{G}(\interior K)&\ge V\left(1-\frac{n}{a_n}-\frac{n}{2}\sum_{i=1}^{n-1}\frac{1}{a_i}\right)\\
&\ge V\left(1-\frac{n}{4}\sum_{i=1}^{n-1}\lambda_i(K_s)-\frac{n}{2}\lambda_n(K_s)\right).
\end{aligned}
$$
The linear normalization in Sect.~\ref{sec:2} preserves interiors and lattice counts, and replaces $V$ by $\vol_n(K)/\det\Lambda$. This proves the result for every full-rank lattice.
\end{proof}

For the standard simplex and an integer $m\ge n+1$, counting positive integer coordinates gives
$$\mathrm{G}(\interior(m\Delta_n))=\binom{m-1}{n}=\frac{m^n}{n!}\left(1-\frac{n(n+1)}{2m}+O(m^{-2})\right)\qquad(m\to\infty),$$
with $n$ fixed. Since every successive minimum is $2/m$, both~\eqref{eq:1-6} and~\eqref{eq:1-7} give the first two terms of this expansion.

\begin{proof}[Proof of Corollary~\ref{cor:1-3}]
Write $\lambda_i=\lambda_i(K_s,\Lambda)$. Since $\lambda_1\le\cdots\le\lambda_n$, pairing $i$ with $n-i$ gives
$$\sum_{i=1}^{n-1}\left(i-\frac{n}{2}\right)\lambda_i=\sum_{1\le i<n/2}\left(\frac{n}{2}-i\right)(\lambda_{n-i}-\lambda_i)\ge0.$$
Thus the subtracted quantity in~\eqref{eq:1-6} is at most $\frac{1}{2}\sum_i i\lambda_i$, proving the corollary. The right-hand sum is empty for $n=1$ and $n=2$.
\end{proof}

\section{The planar product lower bound}\label{sec:8}

We combine two estimates for a fully shifted planar body. Its axial triangle gives a lower bound for the count, while horizontal slices bound the difference between area and count. Retaining the rounding term in the triangle estimate allows these bounds to give the full product.

\begin{lemma}\label{lem:8-1}
Let $p\ge q\ge2$, set $h=\lceil q\rceil-1$ and $\delta=q-h\in(0,1]$, and let $\Delta=\conv(\mathbf{0},p\mathbf{e}_1,q\mathbf{e}_2)$. Then
$$\mathrm{G}(\interior\Delta)\ge\frac{(q-2)(p-\delta)}{2}.$$
\end{lemma}

\begin{proof}
Put $u=p/q-1$, $\varepsilon=1-\delta$, and
$$\rho_k=\lceil(1+u)(k+\delta)\rceil-(1+u)(k+\delta)\qquad(0\le k<h).$$
The positive integer horizontal levels in $\interior\Delta$ are $1,\ldots,h$. Counting their open intervals and reversing the index gives
$$\mathrm{G}(\interior\Delta)=\sum_{k=0}^{h-1}\bigl(\lceil(1+u)(k+\delta)\rceil-1\bigr).$$
Summing the affine terms yields
\begin{equation}\label{eq:8-1}
\mathrm{G}(\interior\Delta)-\frac{(q-2)(p-\delta)}{2}=\sum_{k=0}^{h-1}\rho_k-\frac{h\varepsilon}{2}+\frac{u}{2}\bigl(h+\delta(2-\delta)\bigr).
\end{equation}
If $u\ge\varepsilon$, this is nonnegative because $\rho_k\ge0$.

Suppose $0\le u<\varepsilon$. We show that
\begin{equation}\label{eq:8-2}
\sum_{k=0}^{h-1}\rho_k\ge\frac{h(\varepsilon-u)}{2}.
\end{equation}
For $u=0$, every $\rho_k$ equals $\varepsilon$. For $u>0$, subtracting the integer $k$ from the ceiling argument shows that $\rho_k$ is the ceiling remainder of $\xi_k=\delta+u(k+\delta)$. Since
$$0<\xi_0=(1+u)\delta\le(2-\delta)\delta<1,$$
the initial remainder satisfies $\rho_0=\varepsilon-u\delta\ge\varepsilon-u$. Also $\xi_{k+1}-\xi_k=u\in(0,1)$, so each successive remainder equals either $\rho_k-u$ or $\rho_k+1-u$. Split the sequence at each upward jump. Each resulting run is an arithmetic progression with step $-u$ and nonnegative last term. Its first term is at least $\varepsilon-u$; after a jump it is at least $1-u$. Its mean is therefore at least $(\varepsilon-u)/2$. Summing over the runs proves \eqref{eq:8-2}.

Substitution into \eqref{eq:8-1} gives a lower bound of $u\delta(2-\delta)/2\ge0$, proving the lemma.
\end{proof}

\begin{proof}[Proof of Theorem~\ref{thm:1-4}]
We first prove the required estimate for a fully shifted body $A\subset\R_{\ge0}^2$ with axial intercepts $p\ge q\ge2$. Write
$$A=\{(x,t):0\le t\le q,\ 0\le x\le f(t)\},\qquad V=\vol_2(A),\qquad I=\mathrm{G}(\interior A).$$
Here $f$ is nonnegative and concave, and $f(0)=p$. The case $q=2$ follows from $I\ge0$, so suppose $q>2$. Put $h=\lceil q\rceil-1$ and $\delta=q-h$ as in Lemma~\ref{lem:8-1}. The defining transfer property gives
$$f(s)\ge f(t)+t-s\qquad(0\le s\le t\le q).$$
Hence $g(t)=f(t)-(q-t)$ is nonnegative and nonincreasing. The triangle $\Delta=\conv(\mathbf{0},p\mathbf{e}_1,q\mathbf{e}_2)$ is contained in $A$, so $f(j)>0$ for $1\le j\le h$. The interior count and the integral of $g$ satisfy
$$I=\sum_{j=1}^h(\lceil f(j)\rceil-1)\ge\sum_{j=1}^h(f(j)-1),$$
$$\int_0^q g(t)\,dt\le g(0)+\sum_{j=1}^{h-1}g(j)+\delta g(h)\le g(0)+\sum_{j=1}^h g(j).$$
Substituting $f(t)=q-t+g(t)$ and $h=q-\delta$ yields
$$I\ge V-c_0,\qquad c_0=p+\frac{q}{2}-\frac{\delta(3-\delta)}{2}.$$
On the other hand, Lemma~\ref{lem:8-1} and $\interior\Delta\subseteq\interior A$ give
$$I\ge b_0,\qquad b_0=\frac{(q-2)(p-\delta)}{2}.$$
\begingroup
\postdisplaypenalty=10000
Set $r=(1-1/p)(1-2/q)\in(0,1)$. The two estimates combine to give
$$I\ge(1-r)b_0+r(V-c_0)=rV+\frac{(1-\delta)(q-2)\bigl(q+\delta(p-1)\bigr)}{2pq}\ge rV.$$
The reduction in Sect.~\ref{sec:3} now gives~\eqref{eq:1-8}.
\endgroup
\end{proof}

\begin{remark}\label{rem:8-2}
The stronger symmetric lower bound~\eqref{eq:1-5} in Conjecture~\ref{conj:1-1} fails already for a three-dimensional lattice parallelepiped. Let
$$K=[-2,2]^3,\qquad \Lambda=\{\mathbf{z}\in\Z^3:z_1+z_2+z_3\equiv0\pmod2\}.$$
The cube $[0,2]^3$ with this parity lattice is described in \cite{ref03}; its standard-lattice model is the parallelepiped $M_{12}$ in \cite[Theorem 2.2]{ref04}. The body $K$ is obtained by dilating the cube $[0,2]^3$ by a factor of two and translating it by $-2(1,1,1)\in\Lambda$. Since $\det\Lambda=2$, one has $\vol_3(K)/\det\Lambda=32$. Moreover, all three successive minima are $1/2$. Indeed, $\mathbf{e}_1+\mathbf{e}_2$, $\mathbf{e}_1-\mathbf{e}_2$, and $\mathbf{e}_1+\mathbf{e}_3$ are independent points of $\Lambda\cap\frac{1}{2}K$, while $tK$ contains no nonzero integer point for $t<1/2$. The interior points have coordinates in $\{-1,0,1\}$ and an even number of nonzero coordinates. Hence
$$\mathrm{G}(\interior K,\Lambda)=1+\binom32 2^2=13<\frac{27}{2}=\frac{\vol_3(K)}{\det\Lambda}\prod_{i=1}^3\left(1-\frac{\lambda_i(K,\Lambda)}{2}\right).$$
This disproves the interior-point version~\eqref{eq:1-5} in Conjecture~\ref{conj:1-1}. The original conjecture uses $\mathrm{G}(K,\Lambda)$ \cite[Conjecture 2.2]{ref05}. For $1/2<s<1$ we have $\mathrm{G}(sK,\Lambda)=13$, whereas its conjectured lower bound equals $32(s-1/4)^3>13$ for $s$ sufficiently close to $1$. Thus the original version fails as well.
\end{remark}

\end{document}